\documentclass[11pt,leqno]{article}
\newcommand{\dated}{\mbox{} \hfill {\small [{\tt \today}]}} 
\newtheorem{theorem}{Theorem}[section]
\newtheorem{lemma}[theorem]{Lemma}
\newtheorem{corollary}[theorem]{Corollary}
\newtheorem{proposition}[theorem]{Proposition}
\newtheorem{df}[theorem]{Definition}
\newenvironment{definition}{\begin{df} \rm}{\end{df}}

\usepackage{amsmath,amssymb,amscd}
\newcommand{\pf}[1]{\trivlist \item[\hskip\labelsep\it #1\ ]}
\newcommand{\varpf}[1]{\trivlist \item[\hskip\labelsep\sc #1:]}
\newcommand{\qedbox}{$\rlap{$\sqcap$}\sqcup$}
\newcommand{\qed}{\qquad \qedbox \endtrivlist}
\newcommand{\varqed}{\hfill \rule{0.6em}{0.6em} \endtrivlist}
\newenvironment{proof}{\pf{Proof}}{\qed}

\newenvironment{remark}{\pf{Remark}}{\endtrivlist}
\newenvironment{remarks}{\pf{Remarks} 
   \begin{enumerate}}{\end{enumerate} \endtrivlist}

\newenvironment{items}{
  \begin{enumerate} 
                    
  }{\end{enumerate}}

\newenvironment{keywords}{\noindent\small {\it Keywords\/}:}{\vskip 4pt}
\newenvironment{classification}{\noindent\small 2000 {\it Mathematics Subject
Classification\/}:}{\vskip 12pt}

\renewcommand{\iff}{\quad\Longleftrightarrow\quad}

\newcommand{\comps}{{\mathbb C}}

\newcommand{\posints}{{\mathbb N}}

\newcommand{\void}{\varnothing}

\newcommand{\cstar}{{C^\ast}}

\newcommand{\id}{{\mathrm{id}}}

\newcommand{\A}{{\mathfrak A}}

\begin{document}
\title{Gelfand theory for non-commutative Banach algebras}
\author{{\it Rachid Choukri} \and {\it El Hossein Illoussamen}\thanks{Part of this paper was written while the author was visiting the University of Alberta
in March 2001. Financial support from the University of Alberta is gratefully acknowledged.} \and {\it Volker Runde}\thanks{Research supported by NSERC under grant no.\ 227043-00.}}
\date{}
\maketitle
\begin{abstract}
Let $A$ be a Banach algebra. We call a pair $({\cal G},\A)$ a {\it Gelfand theory\/} for $A$ if the following axioms are
satisfied: (G 1) $\A$ is a $\cstar$-algebra, and ${\cal G} \!: A \to \A$ is a homomorphism; (G 2) the assignment 
$L \mapsto {\cal G}^{-1}(L)$ is a bijection between the sets of maximal modular left ideals of $\A$ and $A$, respectively; (G 3) for each 
maximal modular left ideal $L$ of $\A$, the linear map ${\cal G}_L \!: A / {\cal G}^{-1}(L) \to \A /L $ induced by ${\cal G}$ has dense range.
The Gelfand theory of a commutative Banach algebra is easily seen to be characterized by these axioms. Gelfand theories of arbitrary
Banach algebras enjoy many of the properties of commutative Gelfand theory. We show that
unital, homogeneous Banach algebras always have a Gelfand theory. For liminal $\cstar$-algebras with discrete spectrum,
we show that the identity is the only Gelfand theory (up to an appropriate notion of equivalence).
\end{abstract}
\section{Introduction, motivation, and definition}
Let $A$ be a commutative Banach algebra with character space $\Phi_A$. Then its Gelfand transform ${\cal G}_A$ is the
algebra homomorphism from $A$ into ${\cal C}_0(\Phi_A)$ defined through
\begin{equation} \label{geldef}
  ({\cal G}_A a)(\phi) := \langle a, \phi \rangle \qquad (a \in A, \, \phi \in \Phi_A).
\end{equation}
Without doubt, Gelfand theory is one of the most important tools in the theory of commutative Banach algebras.
\par
At the first glance, Gelfand theory seems to be an exclusively commutative phenomenon: For non-commutative $A$,
the character space $\Phi_A$ is often empty, and even if $\Phi_A \neq \void$, the Gelfand transform
${\cal G}_A \!: A \to {\cal C}_0(\Phi_A)$ as defined in (\ref{geldef}) vanishes on all commutators in $A$,
so that a lot of information may be lost if we pass from $A$ to ${\cal G}A$. It thus doesn't seem to make much
sense to develop Gelfand theory for non-commutative Banach algebras.
\par
The picture changes, however, if we adopt a more abstract point of view: ${\cal C}_0(\Phi_A)$ is a $\cstar$-algebra,
and ${\cal G}_A^\ast \!: {\cal C}_0(\Phi_A)^\ast \to A^\ast$ induces a homeomorphism of $\Phi_{{\cal C}_0(\Phi_A)}$
and $\Phi_A$. The Gelfand transform is thus a homomorphism from $A$ into a $\cstar$-algebra $\A$ which induces a
one-to-one correspondence between the maximal modular ideals of $A$ and $\A$.
\par
For an arbitrary, not necessarily commutative Banach algebra, let $\Lambda_A$ denote the set of maximal modular
left ideals of $A$. For fixed $A$, we consider pairs $({\cal G},\A)$ with the following properties:
\newcounter{gelfand}
\begin{list}{(G \arabic{gelfand})}{\usecounter{gelfand} \settowidth{\labelwidth}{(G 3)}
\settowidth{\leftmargin}{(G 3) }}
\item $\A$ is a $\cstar$-algebra, and ${\cal G} \!: A \to \A$ is a homomorphism.
\item The assignment $\Lambda_\A \ni L \mapsto {\cal G}^{-1}(L)$ is a bijection between $\Lambda_\A$ and
$\Lambda_A$.
\item For each $L \in \Lambda_\A$, the linear map ${\cal G}_L \!: A / {\cal G}^{-1}(L) \to \A /L $ induced by ${\cal G}$ 
has dense range.
\end{list}
\par
If $A$ is commutative, $({\cal G}_A, {\cal C}_0(\Phi_A))$ clearly satisfies (G 1), (G 2), and (G 3). It is less
obvious that these properties already characterize the Gelfand transform of a commutative Banach algebra: 
\begin{proposition} \label{G1}
Let $A$ be a commutative Banach algebra, let $({\cal G}, \A)$ be a pair satisfying\/ {\rm (G 1)}, {\rm (G 2)}, and\/ {\rm (G 3)}.
Then there is an isomorphism $\theta \!: {\mathfrak A} \to {\cal C}_0(\Phi_A)$ such that ${\cal G}_A = \theta \circ {\cal G}$. 
\end{proposition}
\begin{proof}
By (G 2) and (G 3), $\A$ is necessarily commutative. 
Since $\A$ is a commutative $\cstar$-algebra, $\cal G$ is automatically continuous by the classical Gelfand--Rickart theorem.
The adjoint map ${\cal G}^\ast \!: {\mathfrak A}^\ast \to A^\ast$ induces a homeomorphism of the Gelfand spaces 
$\Phi_{\mathfrak A}$ and $\Phi_A$. Identifying $\mathfrak A$ and ${\cal C}_0(\Phi_{\mathfrak A})$ via
the Gelfand--Naimark theorem, we define
\[
  \theta \!: {\mathfrak A} \to {\cal C}_0(\Phi_A), \quad
  f \mapsto f \circ {\cal G}^\ast. 
\]
It is routinely verified that ${\cal G}_A = \theta \circ \cal G$.
\end{proof}
\par
Proposition \ref{G1} motivates the following definition:
\begin{definition} \label{Gdef}
Let $A$ be a Banach algebra.
\begin{items}
\item A pair $({\cal G}, \A)$ satisfying (G 1), (G 2), and (G 3) is called a {\it Gelfand theory\/} for $A$;
the homomorphism ${\cal G}$ is called the corresponding {\it Gelfand transform\/}.
\item We say that two Gelfand theories $({\cal G}_1, \A_1)$ and $({\cal G}_2,\A_2)$ for $A$ are {\it equivalent\/} if there is an isomorphism
$\theta \!: \A_1 \to \A_2$ such that ${\cal G}_2 = \theta \circ {\cal G}_1$.
\item If any two Gefland theories for $A$ are equivalent, we say that $A$ has a {\it unique Gelfand theory\/}.
\end{items}
\end{definition} 
\begin{remarks}
\item By Proposition \ref{G1}, commutative Banach algebras have a unique Gelfand theory, so that we can speak of
{\it the\/} Gelfand transform of a commutative Banach algebra without ambiguity.
\item From (G 2), it is obvious that the kernel of a Gelfand transform is the Jacoboson radical.
\item It does not make much sense in Definition \ref{Gdef}(ii) to require that $\theta$ be a $^\ast$-ho\-mo\-mor\-phism: For example, let 
$A$ be a $\cstar$-algebra, and let $\tau \!: A \to A$ be an automorphism of $A$ which is not a $^\ast$-automorphism. If we had required
that $\theta$ in Definition \ref{Gdef} be an $^\ast$-automorphism, then $(\id_A,A)$ and $(\tau,A)$ would be inequivalent Gelfand theories for $A$.
\end{remarks}
\par
With Definition \ref{Gdef} made, several questions arise naturally: Does every Banach algebra have a Gelfand theory?
If not, which are the Banach algebras with a Gelfand theory? Is a Gelfand theory necessarily unique? Which properties
do Gelfand theories in the sense of Definition \ref{Gdef} have in common with Gelfand theory for commutative Banach
algebras? 
\par
We shall investigate these (and related) questions in this paper.
\subsubsection*{Acknowledgment}
The possibility of a Gelfand theory for homogeneous Banach algebra was indicated to the third author by D.\ W.\ B.\ Somerset during a visit to Aberdeen.
\section{Basic properties of Gelfand transforms}
It is an elementary fact, that the Gelfand transform of a commutative Banach algebra is continuous (in fact: contractive).
In Definition \ref{Gdef}, we have not required a Gelfand transform to be continuous. As our first theorem shows, 
a Gelfand transform is automatically continuous.
\begin{lemma} \label{Gl1}
Let $A$ be a Banach algebra, let $({\cal G}, \A)$ be a Gelfand theory for $A$, let $L \in \Lambda_\A$, and let
$\pi_L \!: \A \to {\cal B}(\A/L)$ the corresponding irreducible representation of $\A$. Then $\pi_L \circ {\cal G}$
is continuous.
\end{lemma}
\begin{proof}
Let $E_L$ again be the image of ${\cal G}_L$ in $\A / L$, so that $E_L$ is a pre-Hilbert space. Since ${\cal G}^{-1}(L)$ is
a maximal modular ideal of $A$ by (G 2), it follows that
\[
  A \to {\cal B}(E_L) , \quad a \mapsto (\pi_L \circ {\cal G})(a)|_{E_L}
\]
is an irreducible representation of $A$ on a normed space and thus continuous by \cite[Theorem 25.7]{BD}, i.e.\
there is $C \geq 0$ such that
\begin{equation} \label{ineq}
  \| (\pi_L \circ {\cal G})(a) |_{E_L} \|_{{\cal B}(E_L)} \leq C \| a \| \qquad (a \in A).
\end{equation}
By (G 3), $E_L$ is dense in $\A / L$ so that the left hand side of (\ref{ineq}) is in fact the operator norm of
$(\pi_L \circ {\cal G})(a)$ on $\A / L$. This yields the claim.
\end{proof}
\begin{theorem} \label{G2}
Let $A$ be a Banach algebra, and let $({\cal G}, \A)$ be a Gelfand theory for $A$. Then ${\cal G}$ is continuous.
\end{theorem}
\begin{proof}
Let $( a_n )_{n=1}^\infty$ be a sequence in $A$ such that there is $b \in \A$ with
\[
  a_n \to 0 \qquad\text{and}\qquad {\cal G}(a_n) \to b.
\]
For $L \in \Lambda_\A$ with corresponding irreducible representation $\pi_L$, the map $\pi_L \circ {\cal G}$ 
is continuous by Lemma \ref{Gl1}, so that
\[
  \pi_L(b) = \lim_{n \to \infty} (\pi_L \circ {\cal G})(a_n) = 0.
\]
Since $\A$ is semisimple, i.e.\ $\bigcap_{L \in \Lambda_\A} \ker \pi_L = \{ 0 \}$, it follows that $b = 0$.
By the closed graph theorem, ${\cal G}$ is continuous.
\end{proof}
\par
If $A$ is a commutative Banach algebra and $a$ is an element of $A$, then the spectra of $a$ in $A$ and
${\cal G}_A a$ in ${\cal C}_0(\Phi_A)$ essentially coincide:
\begin{equation} \label{spectrum}
  \sigma_A(a) \cup \{ 0 \} = \sigma_{{\cal C}_0(\Phi_A)}({\cal G}_A a ) \cup \{ 0 \}.
\end{equation}
Since $\Phi_A$ may be compact for non-unital $A$, (\ref{spectrum}) is the best result we can hope for. It will
turn out that the same is true for Gelfand transfroms in the sense of Definition \ref{Gdef}.
\begin{lemma} \label{Gl2}
Let $A$ be a unital Banach algebra, and let $({\cal G},\A)$ be a Gelfand theory for $A$. Then $\A$ is unital,
and $\cal G$ is a unital homomorphism.
\end{lemma}
\begin{proof}
Let $e \in A$ be the identity of $A$. Let $L \in \Lambda_\A$, and let $\pi_L$ be the corresponding irreducible
represenation of $\A$ on $\A / L$. As in the proof of Lemma \ref{Gl1}, let $E_L$ denote the range of 
${\cal G}_L$ in $\A / L$. It is clear that $(\pi_L \circ {\cal G})(e)$ acts as the identity on $E_L$. Since
$E_L$ is dense in $\A / L$ by (G 3), it follows that $(\pi_L \circ {\cal G})(e) = \id_{\A / L}$. Since
$\A$ is semisimple, this means that ${\cal G}e$ is the identity of $\A$.
\end{proof}
\begin{lemma} \label{Gl3}
Let $A$ be a Banach algebra, and let $({\cal G},\A)$ be a Gelfand theory for $A$. Then, for $a \in A$,
the element ${\cal G}a$ is quasi-invertible in ${\cal G}A$ if and only if it is quasi-invertible in
$\A$.
\end{lemma}
\begin{proof}
Assume that ${\cal G}a$ is quasi-invertible in $\A$, but not in ${\cal G}A$. Consequently, $a$ cannot
be left quasi-invertible in $A$. By \cite[Theorem 2.4.6(e)]{Pal}, this means that there is $\tilde{L} \in \Lambda_A$
such that $a$ is a right modular identity for $\tilde{L}$. By (G 2), there is a unique $L \in \Lambda_\A$ such
that ${\cal G}^{-1}(L) = \tilde{L}$. With $E_L$ and $\pi_L$ as in the proofs of Lemmas \ref{Gl1} and \ref{Gl2},
respectively, let $x := {\cal G}_L(a + \tilde{L}) \in E_L \setminus \{ 0 \}$, and note that
\[
  \pi_L({\cal G}a)x = {\cal G}_L(a^2 + \tilde{L}) = {\cal G}_L(a + \tilde{L}) = x.
\]
Let $b \in \A$ be the quasi-inverse of ${\cal G}a$ in $\A$. Then we have
\[
  x = \pi_L({\cal G}a)x = \pi_L( b {\cal G}a - b)x = \pi_L(b) \pi_L({\cal G}a)x - \pi_L(b)x =
 \pi_L(b)x - \pi_L(b)x =0, 
\]
which is a contradiction.
\end{proof}
\begin{remark}
The proof of Lemma \ref{Gl3}, is essentially a verbatim copy of the proof of \cite[Theorem 4.2.10(c) $\Rightarrow$ (a)]{Pal}.
Note that we cannot just apply \cite[Theorem 4.2.10]{Pal} because ${\cal G}A$ need not be dense in $\A$.
\end{remark}
\par
In perfect analogy with the commutative situation, we have:
\begin{theorem} \label{G3}
Let $A$ be a Banach algebra, and let $({\cal G},\A)$ be a Gelfand theory for $A$. Then, if $A$ is unital,
\begin{equation} \label{spec2}
  \sigma_A(a) = \sigma_\A({\cal G}a) \qquad (a \in A),
\end{equation}
and, if $A$ is non-unital,
\begin{equation} \label{spec3}
  \sigma_A(a) = \sigma_\A({\cal G}a) \cup \{ 0 \} \qquad (a \in A),
\end{equation}
\end{theorem}
\begin{proof}
It is immediate from Lemma \ref{Gl3} that
\[
  \sigma_A(a) \cup \{ 0 \} = \sigma_\A({\cal G}a) \cup \{ 0 \} \qquad (a \in A),
\]
no matter if $A$ has an identity or not. If $A$ is non-unital, $0 \in \sigma_A(a)$ for every $a \in A$. This
establishes (\ref{spec3}).
\par
Suppose that $A$ is unital. If $0 \notin \sigma_A(a)$, i.e.\ if $a$ is invertible in $A$, then ${\cal G}a$ is invertible in
$\A$ by Lemma \ref{Gl2}, so that $0 \notin \sigma_\A({\cal G}a)$. Let $a \in A$ be such that $0 \notin \sigma_\A({\cal G}a)$,
and assume that $0 \in \sigma_A(a)$, i.e.\ $a$ is not invertible in $A$. We first treat the case where $a$ is not left
invertible. Then there is $\tilde{L} \in \Lambda_A$ with $a \in L$, and by (G 2), there is a unique $L \in \Lambda_\A$
such that $\tilde{L} = {\cal G}^{-1}(L)$. This, however, means that $a \in L$, which is impossible because ${\cal G}a$
is invertible in $\A$. Assume now that $a$ has a left inverse $b \in A$. Then ${\cal G}b = ({\cal G}a)^{-1}$, so that
$0 \notin \sigma_\A({\cal G}b)$. Since $b$ has $a$ as its right inverse, $b$ cannot be left invertible in $A$ 
(otherwise $b$ would be invertible with inverse $a$ and thus $0 \notin \sigma_A(a)$). But as we have just seen, this
is impossible.
\end{proof}
\begin{remark}
In both Theorem \ref{G2} and \ref{G3}, we have not made use of the fact that $\A$ is a $\cstar$-algebra. All
we need is that $\A$ is semisimple.
\end{remark}
\section{Existence of Gelfand theories}
Which Banach algebras do have Gelfand theories? Of course, the commutative ones. Also, if $A$ is a $\cstar$-algebra, then
$(\id_A,A)$ is trivially a Gelfand theory for $A$. What about non-trivial examples? 
\par
Our first result is negative; it shows that there are Banach algebras without Gelfand theories.
\begin{proposition} \label{G4}
Let $E$ be the Banach space $c_0$ or $\ell^p$ with $p \in (1,\infty) \setminus \{ 2 \}$. Then there is no Gelfand theory 
for ${\cal B}(E)$.
\end{proposition}
\begin{proof}
Assume that there is a $\cstar$-algebra $\mathfrak A$ and a homomorphism ${\cal G} \!: {\cal B}(E) \to \mathfrak A$ such that 
(G 2) and (G 3) are satisfied. Fix $x \in E \setminus \{ 0 \}$. Then $\tilde{L} := \{ T \in {\cal B}(E) : T x = 0 \}$ 
is a maximal 
modular left ideal of ${\cal B}(E)$; note that $A / \tilde{L} \cong E$. By assumption, there is a maximal modular left ideal 
$L$ of $\mathfrak A$ such that $\tilde{L} = {\cal G}^{-1}(L)$. Since, by Theorem \ref{G2}, ${\cal G}$ is continuous, the same 
is true for ${\cal G}_L \!: E \cong A /\tilde{L} \to \A / L$. Hence, since $E$ is separable, the Hilbert space $\A/L$ must be 
separable as well, i.e.\ $\A / L \cong \ell^2$. Therefore, $\cal G$ induces a non-zero homomorphism from ${\cal B}(E)$ into 
${\cal B}(\ell^2)$. This is impossible by \cite[6.15 Corollary]{BP}.
\end{proof}
\par
Our first positive result is about homogeneous Banach algebras. Recall that a Banach algebra $A$ is called 
{\it $n$-homogeneous\/} for some $n \in \posints$ if all its irreducible representations have dimension $n$. 
If we don't want to emphasize for which $n \in \posints$ the algebra $A$ is $n$-homogeneous, we call $A$ simply 
{\it homogeneous\/}. 
\begin{lemma} \label{homlem3}
Let $n \in \posints$, let $A$ be an $n$-homomogeneous Banach algebra, and let $P \in \Pi_A$, the set of all primitive ideals of $A$. Then there is an algebra homomorphism 
$\pi_P \!: A \to {\mathbb M}_n$ with $P = \ker \pi_P$ and $\| \pi_P \| \leq \sqrt{n}$.
\end{lemma}
\begin{proof}
Let $\pi \!: A \to A/P$ be the quotient map. Since $A$ is $n$-homogeneous, we have (algebraically) $A/P \cong {\mathbb M}_n$. By \cite[Lemma 2]{DR}, there is an isomorphism
$\alpha \!: A/P \to {\mathbb M}_n$ with $\| \alpha \| \leq \sqrt{n}$. Let $\pi_P := \alpha \circ \pi$.
\end{proof}
\begin{theorem} \label{hom}
Let $A$ be a unital homogeneous Banach algebra. Then there is a Gelfand theory for $A$.
\end{theorem}
\begin{proof}
Let $n \in \posints$ be such that $A$ is $n$-homogeneous. For each $P \in \Pi_A$, let $\pi_P \!: A \to {\mathbb M}_n$ be as in Lemma \ref{homlem3}.
Define
\[
  {\cal G} \!: A \to \text{$\ell^\infty$-}\bigoplus_{P \in \Pi_A} {\mathbb M}_n, \qquad a \mapsto (\pi_P(a))_{P \in \Pi_A},
\]
and let $\A$ denote the $\cstar$-subalgebra of $\text{$\ell^\infty$-}\bigoplus_{P \in \Pi_A} {\mathbb M}_n$ generated by ${\cal G}A$. Clearly, $({\cal G},\A)$
satisfies (G 1).
\par
Let $P \in \Pi_\A$. Since $\A$ --- as a subalgebra of $\text{$\ell^\infty$-}\bigoplus_{P \in \Pi_A} {\mathbb M}_n$ --- satisfies the polynomial identity $S_{2n} = 0$,
it follows from \cite[Lemma 3.4]{AD} that $A/P \cong {\mathbb M}_m$ for some $m \in \posints$ with $m \leq n$. Let $\pi \!: \A \to \A / P$ be the quotient map. Then
$\pi \circ {\cal G}$ is a non-zero algebra homomorphism from $A$ into $A/P \cong {\mathbb M}_m$. Since $A$ is $n$-homomogeneous, and since $m \leq n$, it follows that
${\cal G}^{-1}(P) = \ker(\pi \circ {\cal G}) \in \Pi_A$ (so that, in fact, $m = n$). We thus have a map
\begin{equation} \label{primmap}
  \Pi_\A \to \Pi_A, \quad P \mapsto {\cal G}^{-1}(P).
\end{equation}
From the definition of $A$, it is clear that (\ref{primmap}) is surjective. Let $P_1, P_2 \in \Pi_\A$ with $A \cap {\cal G}^{-1}(P_j) = A \cap {\cal G}^{-1}(P_2)$, and let 
$\pi_j \!: \A \to {\mathbb M}_n$ be $^\ast$-homomorphisms such that $P_j = \ker \pi_j$ for $j=1,2$. The each $\pi_j$ induces a unital contractive homomorphism from
$\pi_P(A)$ to ${\mathbb M}_n$. As unit preserving contractions, these homomorphisms are necessarily positive and thus $^\ast$-homomorphisms. Hence, they are unitarily equivalent. It follows that $\pi_1$ and $\pi_2$
are unitarily equivalent, so that $P_1 = \ker \pi_1 = \ker \pi_2 = P_2$. Consequently, (\ref{primmap}) is also injective. From the preceeding discussion, it has become clear
that $\A$ is also $n$-homogeneous. Hence, for each $P \in \Pi_\A$, we have a canonical isomorphism
\[
  A/ {\cal G}^{-1}(P) \cong \A / P.
\]
It follows that $({\cal G},\A)$ satisfies (G 2) and (G 3) as well and thus is a Gelfand theory for $A$.
\end{proof}
\par
Let $A$ be a Banach $^\ast$-algebra. Then there is a largest $C^\ast$-seminorm $\gamma_A$ on $A$; the completion of
$A / \ker \gamma_A$ with respect to $\gamma_A$ is called the {\it enveloping $\cstar$-algebra\/} of $A$ and is denoted
by $\cstar(A)$. If, for example, $A = L^1(G)$ for some locally compact group $G$, the enveloping $\cstar$-algebra
$A$ is just the usual group $\cstar$-algebra. We denote the canonical $^\ast$-homomorphism from $A$ into $\cstar(A)$
by $\iota_A$. Recall that $A$ is said to be {\it hermitian\/} if each self-adjoint element of $a$ has real spectrum in
$A$. It is thus immediate from Theorem \ref{G3} that $(\iota_A,\cstar(A))$ is a Gelfand theory for $A$ only if
$A$ is hermitian. 
\par
In fact, the converse holds as well:
\begin{proposition}
Let $A$ be a Banach $^\ast$-algebra. Then $(\iota_A, \cstar(A))$ is a Gelfand theory for $A$ if and only if $A$ is 
hermitian.
\end{proposition}
\begin{proof}
Suppose that $A$ is hermitian. Let $L \in \Lambda_{\cstar(A)}$, and let $\phi$ be a pure state on $\cstar(A)$ 
such that
\[
  L = \{ a \in \cstar(A) : \langle a^\ast a , \phi \rangle = 0 \}
\]
(see \cite[Theorem 5.3.5]{Mur}). By \cite{Pal}, we have $\iota_A^{-1}(L) \in \Lambda_A$, and every maximal modular ideal of $A$ arises in this fashion. Since
$\iota_A(A)$ is dense in $\cstar(A)$, a pure state on $\cstar(A)$ is already determined by the values it takes on $A$.
This establishes (G 2). The density of $\iota_A(A)$ in $\cstar(A)$ also implies (G 3).
\end{proof}
\section{Hereditary properties}
What happens to the existence of Gelfand theories when we pass from a Banach algebra to a quotient, to a closed ideal, etc.?
\par
We first deal with the question of how the existence of Gelfand theories is affected if we adjoin an identity. The following algebraic lemma
is standard, but for the reader's convenience, we outline a proof:
\begin{lemma}
Let $A$ be an algebra, let $I$ be an ideal of $A$, and let $\pi$ be an irreducible representation of $I$ on a linear space $E$. Then $\pi$ extends to a unique
irreducible representation of $A$ on $E$. Conversely, if $\pi$ is an irreducible representation of $A$ on a linear space $E$, then $\pi |_I$ is an irreducible
representation of $I$ on $E$.
\end{lemma}
\begin{proof}
Let $a \in A$, and let $x \in E$. Since $\pi$ is irreducible, there are $b \in I$ and $y \in E$ such that $x = \pi(b)y$. Define $\pi(a) x := \pi(ab) x$. It is routinely checked
that this definition is independent of the choices of $b$ and $x$ and yields an irreducible representation of $A$ on $E$. Assume that $\pi$ has two such extensions $\pi_1$ and
$\pi_2$. Let $a \in A$, and $x = \pi(b)y$ with $b \in I$ and $y \in E$. Then we have
\[
  \pi_1(a) x = \pi_1(a)\pi(b)y = \pi(ab)y = \pi_2(a)\pi(b) y = \pi_2(a) x,
\]
so that $\pi_1 = \pi_2$.
\par
Let $\pi$ be an irreducible representation of $A$ on a linear space $E$. Let $x \in E$. Since $\pi(I) x$ is invariant under $\pi(A)$, it follows that $\pi(I) x = E$ or
$\pi(I)x = \{ 0 \}$. Hence, $\pi |_I$ is irreducible.
\end{proof}
\par
The correspondence between irreducible representations and maximal modular left ideals (\cite[Proposition 24.5]{BD}) then yields immediately:
\begin{corollary} \label{idcor}
Let $A$ be a Banach algebra, and let $I$ be a closed ideal of $A$. Then
\[
  \{ L \in \Lambda_A : I \not\subset L \} \to \Lambda_I, \quad L \mapsto I \cap L
\]
is a bijection.
\end{corollary}
\par
We apply Corollary \ref{idcor} to establish first hereditary properties of the existence of Gelfand theories:
\begin{proposition} \label{unitprop}
Let $A$ be a non-unital Banach algebra which has a Gelfand theory. Then $A^\#$ has also a Gelfand theory.
\end{proposition}
\begin{proof}
Let $({\cal G},\A)$ be a Gelfand theory for $A$. Let $\A^\#$ denote the {\it unconditional\/} unitization of $\A$, i.e.\ if $\A$ already has an identity, we adjoin another one.
Define a homomorphism
\[
  {\cal G}^\# \!: \A^\# \to \A^\#, \quad a + \lambda e_{A^\#} \mapsto {\cal G}a + \lambda e_{\A^\#}.
\] 
By Corollary \ref{idcor}, $({\cal G}^\#, \A^\#)$ satisfies (G 2).
\par
If $\A$ is non-unital, one obtains a $\cstar$-norm on $\A^\#$ through a faithful representation of $\A$ on some Hilbert space. If $\A$ is unital, we have a $^\ast$-isomorphism
\[
  \A^\# \to \A \oplus \comps, \quad a + \lambda e_{\A^\#} \mapsto a + \lambda e_\A \oplus \lambda, 
\]
which again endows $\A^\#$ with a $\cstar$-norm. Hence, $({\cal G}^\#,\A^\#)$ satisfies (G 1) as well.
\par
It is routinely verified that (G 3) also holds, so that $({\cal G}^\#,\A^\#)$ is indeed a Gelfand theory for $A^\#$.
\end{proof}
\begin{proposition} \label{idprop}
Let $A$ be a Banach algebra which has a Gelfand theory, and let $I$ be a closed ideal of $A$. Then $I$ has a Gelfand theory.
\end{proposition}
\begin{proof}
Let $({\cal G},\A)$ be a Gelfand theory for $A$. Let
\[
  {\mathfrak I} := \bigcap \{ L \in \Lambda_\A : {\cal G}I \subset L \}.
\]
Obviously, $\mathfrak I$ is a closed left ideal of $\A$. We claim that $\mathfrak I$ is in fact a two-sided ideal of $\A$ and thus a $\cstar$-algebra.
Let $L \in \Lambda_\A$ be such that ${\cal G}I \subset L$, and let
\[
  P := \{ a \in \A : a \A \subset L \}.
\]
It is clear that
\begin{equation} \label{incl}
  {\cal G}^{-1}(P) \subset \{ a \in A : a A \subset {\cal G}^{-1}(L) \} =: Q.
\end{equation}
Conversely, let $a \in Q$, and let $\pi_L$ denote the irreducible representation of $\A$ on $\A / L$. Let $E_L$ be the image of ${\cal G}_L$ in $\A /L$. Then
$\pi_L({\cal G} a)$ vanishes on $E_L$ by the definition of $Q$. Since $E_L$ is dense in $\A / L$ by (G 3), it follows that $\pi_L({\cal G}a) = 0$, i.e.\ $a \in {\cal G}^{-1}(P)$.
Hence, the inclusion (\ref{incl}) is an equality.
Since $I \subset {\cal G}^{-1}(L)$, it follows that $I \subset {\cal G}^{-1}(P)$ and thus ${\cal G}I \subset P$. Hence, 
\[
  {\mathfrak I} = \bigcap \{ P \in \Pi_\A : {\cal G}I \subset P \}
\]
is a two-sided ideal of $\A$.
\par
From Corollary \ref{idcor}, it is now straightforward that $({\cal G} |_I, {\mathfrak I})$ satisfies (G 2).
Furthermore, for any $L \in \Lambda_\A$, we have canonical isomorphisms
\[
  {\mathfrak I} / {\mathfrak I} \cap L \cong \A / L \qquad\text{and}\qquad I / I \cap {\cal G}^{-1}(L) \cong A / {\cal G}^{-1}(L).
\]
This establishes (G 3) for $({\cal G}|_I, {\mathfrak I})$.
\end{proof}
\par
Combining Propositions \ref{unitprop} and \ref{idprop}, we obtain:
\begin{corollary} \label{unitcor}
Let $A$ be a non-unital Banach algebra. Then $A$ has a Gelfand theory if and only if $A^\#$ has one.
\end{corollary}
\par
We next turn to quotients:
\begin{proposition}
Let $A$ be a Banach algebra with a Gelfand theory, and let $I$ be a closed ideal of $A$. Then $A / I$ has a Gelfand theory.
\end{proposition}
\begin{proof}
Let $({\cal G},\A)$ be a Gelfand theory for $A$, and let again
\[
  {\mathfrak I} := \bigcap \{ L \in \Lambda_\A : {\cal G}I \subset L \}.
\]
Let $\pi \!: \A \to \A / {\mathfrak I}$ be the quotient map. By the definition of $\mathfrak I$, $\pi \circ {\cal G}$ vanishes on $I$ and thus induces
an algebra homomorphism $\tilde{G} \!: A/I \to \A / {\mathfrak I}$. With Corollary \ref{idcor}, it is easy to see that $(\tilde{\cal G}, \A / {\mathfrak I})$
is a Gelfand theory for $A/I$.
\end{proof}
\section{Gelfand theory for $\cstar$-algebras}
If $A$ is a $\cstar$-algebra, then $(\id_A,A)$ is trivially a Gelfand theory for $A$. But what can we say about arbitrary Gelfand theories $({\cal G},\A)$ for $A$?
Is a Gelfand theory necessarily equivalent to $(\id_A,A)$? What if we impose additional conditions on $\cal G$ and $\A$?
\par
Let $A$ be a $\cstar$-algebra, and let $({\cal G},\A)$ be a Gelfand theory for $A$. Then, for each $L \in \Lambda_\A$, the spaces $A / {\cal G}^{-1}(L)$ and $\A / L$
are Hilbert spaces (\cite[2.8.5 Corollary]{Dix}). For any Hilbert space $\mathfrak H$ and $\xi, \eta \in {\mathfrak H}$, define a rank one operator
\[
  \xi \odot \eta \!: {\mathfrak H} \to {\mathfrak H}, \quad x \mapsto \langle x, \eta \rangle \xi.
\]
\begin{lemma} \label{rankone}
Let $A$ be a $\cstar$-algebra, let $({\cal G}, \A)$ be a Gelfand theory for $A$, and let $a \in A$ be such that $\pi_{{\cal G}^{-1}(L)}(a) = \xi \odot {\cal G}_L^\ast \eta$ with
$\xi \in A / {\cal G}^{-1}(L)$ and $\eta \in \A / L$. Then $\pi_L({\cal G}a) = {\cal G}_L \xi \odot \eta$ holds.
\end{lemma}
\begin{proof}
Since
\begin{equation} \label{intertw}
  {\cal G}_L \circ \pi_{{\cal G}^{-1}(L)}(a) = \pi_L({\cal G}a) \circ {\cal G}_L,
\end{equation}
and since ${\cal G}_L$ has dense range, it follows that $\pi_L({\cal G}a)$ is also a rank one operator, whose range is spanned by ${\cal G}_L \xi$. Hence, there is a unique
$\tilde{\eta} \in \A / L$ such that $\pi_L({\cal G}a) = {\cal G}_L \xi \odot \tilde{\eta}$. By (\ref{intertw}), we have
\[
  \langle x, {\cal G}_L^\ast \eta \rangle = \langle {\cal G}_L x, \eta \rangle \qquad (a \in A / {\cal G}^{-1}(L))
\]
and thus $\tilde{\eta} = \eta$.
\end{proof}
\par
Recall that a $\cstar$-algebra $A$ is called {\it liminal\/} if $\pi(A) = {\cal K}({\mathfrak H})$ for each irreducible $^\ast$-representation of $A$ on Hilbert space $\mathfrak H$. If
we only have $\pi(A) \supset {\cal K}({\mathfrak H})$, then we say that $A$ is {\it postliminal\/} (see \cite{Dix}).
\begin{corollary} \label{limcor}
Let $A$ be a $\cstar$-algebra, and let $({\cal G},\A)$ be a Gelfand theory for $A$. Then, if $A$ is postliminal, so is $\A$, and if $A$ is liminal, so is $\A$.
\end{corollary}
\begin{proof}
Suppose that $A$ is postliminal, and let $L \in \Lambda_\A$. By Lemma \ref{rankone}, we thus have $\pi_L(\A) \cap {\cal K}(\A/L) \neq \{ 0 \}$. Since $\pi_L$ is irreducible, this
already implies that ${\cal K}(\A/L) \subset \pi_L(\A)$.
\par
Suppose that $A$ is liminal. By the first part of the proof, $\A$ is postliminal. Let $\mathfrak I$ be the largest, closed, liminal ideal of $\A$, and assume that ${\mathfrak I} \subsetneq \A$.
Then there is $L \in \Lambda_\A$ such that ${\mathfrak I} \subset L$. Since ${\cal G}_L$ is injective, so is ${\cal G}_L^\ast$. It follows that the closed linear span of
$\{ \xi \odot {\cal G}_L^\ast \eta : \xi \in A / {\cal G}^{-1}(L), \, \eta \in \A / L \}$ is all of ${\cal K}(A / {\cal G}^{-1}(L))$. From Lemma \ref{rankone} again, we conclude that
${\cal G}A \subset {\mathfrak I}$ and thus ${\cal G}^{-1}(L) = A \notin \Lambda_A$. This contradicts (G 2).
\end{proof}
\par
This enables us to prove a first, modest uniqueness result for Gelfand theories:
\begin{proposition}
Let $A$ be a unital, postliminal $\cstar$-algebra, and let $({\cal G},\A)$ be a Gelfand theory for $A$ such that $\| {\cal G} \| \leq 1$. Then $({\cal G},\A)$ and
$(\id_A,A)$ are equivalent.
\end{proposition}
\begin{proof}
Since $A$ is unital, so is $\A$, and $\cal G$ is a unital homomorphism. Hence, $\cal G$ is an injective $^\ast$-homomorphism, and ${\cal G}A$ is a $\cstar$-subalgebra of
$\A$. By Corollary \ref{limcor}, $\A$ is also postliminal. From (G 2), it follows that the hypotheses of the non-commutative Stone--Weierstra{\ss} theorem \cite[11.1.8 Th\'eor\`eme]{Dix}
are satisfied. Hence, we have $\A = {\cal G}A$.
\end{proof}
\par
To obtain better uniqueness results, we require both stronger tools and further restrictions on the $\cstar$-algebra $A$.
\par
Our next lemma is a subscript to a deep theorem by U.\ Haagerup (\cite[Theorem 7.5]{Pis}):
\begin{lemma} \label{cyc}
Let $A$ be a $\cstar$-algebra, and let $\pi$ be a bounded representation of $A$ on a $\mathfrak H$ such that there is a cyclic vector $\xi \in \mathfrak H$ for $\pi$, i.e.\
$\pi(A) \xi$ is dense in $\mathfrak H$. Then there is an invertible operator $S \in {\cal B}({\mathfrak H})$ with $\| S \| \| S^{-1} \| \leq (1+2 \| \pi \|)^4$ such that
\[
  A \to {\cal B}({\mathfrak H}), \quad a \mapsto S \pi(a) S^{-1}
\]
is a $^\ast$-representation.
\end{lemma}
\begin{proof}
If $A$ is unital, the presence of a cyclic vector forces $\pi$ to be unital. Hence, in this case, the claim follows immediately from \cite[Theorem 7.5]{Pis} (with a much better norm estimate).
\par
Suppose that $A$ is not unital. Adjoin an identity, and let $\pi^\# \!: A^\ast \to {\cal B}({\mathfrak H})$ be the unital extension of $\pi$. Let $\phi \in \Phi_{A^\#}$ be the character 
corresponding to the maximal ideal $A$. Then we have
\[
  \| \pi^\#(a) \| = \| \langle a,\phi \rangle \id_{\mathfrak H} + \pi(a - \langle a,\phi \rangle) \| \leq (1 + 2 \| \pi \|) \| a \| \qquad (a \in A^\#).
\]
The claim then follows again from \cite[Theorem 7.5]{Pis}.
\end{proof}
\begin{lemma} \label{limiso}
Let $A$ be a liminal $\cstar$-algebra, and let $({\cal G},\A)$ be a Gelfand theory for $A$. Then, for each $L \in \Lambda_\A$, the Gelfand transform $\cal G$ induces
an isomorphism of ${\cal K}(A /{\cal G}^{-1}(L))$ and ${\cal K}(\A/L)$ whose inverse has norm at most $(1 + 2 \| {\cal G} \|)^4$.
\end{lemma}
\begin{proof}
As we have already remarked in the proof of Corollary \ref{limcor}, the closed linear span of $\{ \xi \odot {\cal G}_L^\ast \eta : \xi \in A / {\cal G}^{-1}(L), \, \eta \in \A / L \}$
is dense in ${\cal K}(A/ {\cal G}^{-1}(L))$, it follows that the homomomorphism $\tilde{\cal G} \!: {\cal K}(A/ {\cal G}^{-1}(L)) \to {\cal K}(\A/L)$ induced by $\cal G$ has dense range.
By (G 3) every non-zero vector in the range of ${\cal G}_L$ is cyclic for $\tilde{G}$. By Lemma \ref{cyc}, this means that there is an invertible operator $S \in {\cal B}(\A/L)$ with
$\| S \| \| S^{-1} \| \leq (1 + 2 \| \tilde{\cal G} \|)^4 \leq (1 + 2 \| {\cal G} \|)^4$ such that
\begin{equation} \label{asthom}
  {\cal K}(A/ {\cal G}^{-1}(L)) \to {\cal K}(\A/L), \quad a \mapsto S (\tilde{\cal G}a) S^{-1}
\end{equation}
is a $^\ast$-homomorphism and thus has closed range. Hence, (\ref{asthom}) and $\tilde{G}$ are in fact isomorphisms. The claim about $\| \tilde{\cal G}^{-1} \|$ follows from the estimate
for $\| S \| \| S^{-1} \|$ and the fact the $^\ast$-isomorphisms are always isometric.
\end{proof}
\par
For the next lemma, recall that, for any Banach algebra $A$, the set $\Pi_A$ of primitive ideals of $A$ is canonically equipped with its {\it Jacoboson topology\/}
(\cite[Definition 26.3]{BD}).
\begin{lemma} \label{primlem}
Let $A$ be a liminal $\cstar$-algebra, and let $({\cal G},\A)$ be a Gelfand theory for $A$. Then
\begin{equation} \label{primitive}
  \Pi_\A \to \Pi_A, \quad P \mapsto {\cal G}^{-1}(P)
\end{equation}
is a continuous bijection.
\end{lemma}
\begin{proof}
It was shown implicitly in the proof of Proposition \ref{idprop} that (\ref{primitive}) is well-defined and surjective (this is true for any Banach algebra $A$).
\par
Let $P_1, P_2 \in \Pi_\A$ be such that ${\cal G}^{-1}(P_1) = {\cal G}^{-1}(P_2) =: P$, and assume that $P_1 \neq P_2$. By Lemma \ref{limiso}, $\cal G$ induces isomorphism 
${\cal G}_j \!: A / P \to \A/ P_j$ for $j =1,2$. Let $L \in \Lambda_{A/P} \subset \Lambda_A$ (with the appropriate identification). Let $\pi \!: \A \to \A /P_1 \cap P_2 \cong \A/P_1 \oplus \A/P_2$.
Let $L_1 := \pi^{-1}(\A / P_1 \oplus {\cal G}_2(L))$ and $L_2 := \pi^{-1}({\cal G}(L_1) \oplus \A /P_2)$. It follows that $L = {\cal G}^{-1}(L_j)$ for $j=1,2$, which violates the injectivity
hypothesis in (G 2). Hence, (\ref{primitive}) is injective.
\par
Let $F \subset \Pi_A$ be closed. Then there is a closed ideal $I$ of $A$ such that $F = \{ P \in \Pi_A : I \subset P \}$. Let
\[
  {\mathfrak I} := \bigcap \{ L \in \Lambda_\A : {\cal G}I \subset L \}.
\]
In the proof of Proposition \ref{idprop}, we have seen that $\mathfrak I$ is an ideal of $\A$. Let $P \in \Lambda_\A$. By Lemma \ref{limiso}, we have
\begin{eqnarray*}
  {\cal G}^{-1}(P) = \bigcap \{ L \in \Lambda_A : {\cal G}^{-1} \subset L \} = \bigcap \{ {\cal G}^{-1}(L) : L \in \Lambda_\A, \, P \subset L \}.
\end{eqnarray*}
This implies for $P \in \Pi_\A$ that
\[
  {\cal G}^{-1}(P) \in F \iff P \supset {\mathfrak I}.
\]
Hence, the inverse image of $F$ under (\ref{primitive}) is $\{ P \in \Pi_\A : {\mathfrak I} \subset P \}$ and thus closed in $\Pi_A$.
\end{proof}
\begin{theorem} \label{limthm}
Let $A$ be a liminal $\cstar$-algebra such that $\Pi_A$ is discrete. Then $A$ has a unique Gelfand theory.
\end{theorem}
\begin{proof}
Let $({\cal G},\A)$ be a Gelfand theory for $A$. By Corollary \ref{limcor}, $\A$ is liminal as well. Let $P \in \Pi_\A$. Since $\Pi_A$ is discrete, $\{ {\cal G}^{-1}(P) \}$ is open in 
$\Pi_A$. The continuity assertion of Lemma \ref{primlem}, then implies that $\{ P \}$ is open in $\Pi_\A$. Hence, $\Pi_\A$ is discrete as well.
\par
The general theory of liminal $\cstar$-algebras with Hausdorff spectrum (\cite[Chapter 10]{Dix}) implies that there are families $({\mathfrak H}_\alpha )_\alpha$ and
$( {\mathfrak K}_\alpha )_\alpha$ of Hilbert spaces such that
\[
  A \cong \text{$c_0$-} \bigoplus_\alpha {\cal K}({\mathfrak H}_\alpha) \qquad\text{and}\qquad  A \cong \text{$c_0$-} \bigoplus_\alpha {\cal K}({\mathfrak K}_\alpha).
\]
By Lemma \ref{limiso}, $\cal G$ induces, for each index $\alpha$ an isomorphism ${\cal G}_\alpha \!: {\cal K}({\mathfrak H}_\alpha) \to {\cal K}({\mathfrak K}_\alpha)$ such 
that $\| {\cal G}_\alpha^{-1} \| \leq (1 + 2 \| {\cal G} \|)^4$. It follows that $\cal G$ is bounded below and thus an isomorphism.
\end{proof}
\par
Without any restrictions placed on the $\cstar$-algebra under consideration, the following is the strongest assertion we can make about Gelfand theories of $\cstar$-algebras:
\begin{proposition} \label{cstarprop}
Let $A$ be a $\cstar$-algebra and let $({\cal G},\A)$ be a Gelfand theory for $A$. Then ${\cal G}A$ is closed in $\A$.
\end{proposition}
\begin{proof}
For each $L \in \Lambda_\A$, let $\pi_L$ denote the corresponding irreducible $^\ast$-representation of $\A$ on $\A /L$, and let $E_L$
be the image of ${\cal G}_L$ in $\A /L$. Define $\rho_L := \pi_L \circ {\cal G}$. Then $\rho_L$ is a bounded representation of $A$ on $\A/L$ with $\| \rho_L \| \leq
\| {\cal G} \|$. As in the proof of Lemma \ref{limiso}, we obtain an invertible operator $S_L \in {\cal B}(\A/L)$ with $\| S_L \| \| S_L^{-1} \| \leq (1+2\| {\cal G} \|)^4$ such that
\[
  A \to {\cal B}(\A/L), \quad a \mapsto S_L \rho_L(a) S_L^{-1}
\]
is a $^\ast$-homomorphism. Define
\[
  \theta \!: \A \to \text{$\ell^\infty$-}\bigoplus_{L \in \Lambda_\A} {\cal B}(\A/L), \qquad a \mapsto (S_L \pi_L(a) S_L^{-1})_{L \in \Lambda_\A}.
\]
Then $\theta \circ {\cal G} \!: A \to \text{$\ell^\infty$-}\bigoplus_{L \in \Lambda_\A} {\cal B}(\A/L)$ is a $^\ast$-homomorphism, whose range is a $\cstar$-subalgebra of 
$\text{$\ell^\infty$-}\bigoplus_{L \in \Lambda_\A} {\cal B}(\A/L)$, which we may identify with $A$.  
\par
Let $( a_n )_{n=1}^\infty$ be a sequence in $A$ such that ${\cal G}a_n \to b \in \A$. Since $a_n = \theta({\cal G}a_n)$ for $n \in \posints$, the sequence $( a_n )_{n=1}^\infty$ is 
Cauchy in $A$ and thus convergent to some $a \in A$. It is clear that $b = {\cal G}a$. 
\end{proof}
\section{Outlook}
The results we have obtained in this paper certainly leave room for improvement. Here are a few questions which are natural, but which we have been unable to settle:
\begin{enumerate}
\item Let $({\cal G},\A)$ be a Gelfand theory for a Banach algebra $A$. Does ${\cal G}A$ generate $\A$ as a $\cstar$-algebra (at least if $\A$ is postliminal)?
\item Let $A$ be a Banach algebra, and let $I$ be a closed ideal of $A$ such that $I$ and $A/I$ both have Gelfand theories. Does $A$ have a Gelfand theory?
\item Does a homogeneous Banach algebra have a unique Gelfand theory?
\item Does a $\cstar$-algebra $A$ have a unique Gelfand theory (at least if $A$ is postliminal/liminal/liminal with Hausdorff spectrum)? 
\end{enumerate}
\par
Also, the axioms (G 1), (G 2), and (G 3) were chosen mainly to allow for a simultaneous treatment of commutative Gelfand theory and the enveloping $\cstar$-algebra of
a Banach $^\ast$-algebra. Possibly, in order to develop a farther reaching non-commutative Gelfand theory, one has to consider a different set of axioms.
\dated
\vfill
\noindent
\begin{tabbing}
{\it Second author's address\/}: \= Department of Mathematical Sciences \kill 
{\it First author's address\/}:  \> D\a'epartement de Math\a'ematiques \\
                                 \> Ecole Normale Sup\a'erieure Takaddoum \\
                                 \> B.P.\ 5118 \\
                                 \> 10105 Rabat \\[\smallskipamount]
{\it E-mail\/}:                  \> {\tt rchoukri@hotmail.com} \\[\bigskipamount]
{\it Second author's address\/}: \> D\a'epartement de Math\a'ematiques \\
                                 \> Ecole Normale Sup\a'erieure Takaddoum \\
                                 \> B.P.\ 5118 \\
                                 \> 10105 Rabat \\[\smallskipamount]
{\it E-mail\/}:                  \> {\tt illous@hotmail.com} \\[\bigskipamount]
{\it Third author's address\/}:  \> Department of Mathematical Sciences \\
                                 \> University of Alberta \\
                                 \> Edmonton, Alberta \\
                                 \> Canada T6G 2G1 \\[\smallskipamount]
{\it E-mail\/}:                  \> {\tt runde@math.ualberta.ca} \\
                                 \> {\tt vrunde@ualberta.ca} \\[\smallskipamount]
{\it URL\/}:                     \> {\tt http://www.math.ualberta.ca/$^\sim$runde/runde.html} \
\end{tabbing}
\end{document}